\documentclass[fleqn,final,3p,times,authoryear]{elsarticle}
\usepackage{float}
\usepackage{graphicx}
\usepackage{amsmath}
\usepackage{amssymb}
\usepackage{amsthm}
\usepackage{mathtools}
\usepackage{cancel}
\usepackage{enumerate}
\usepackage{fourier}
\usepackage[title]{appendix}
\usepackage[usenames,dvipsnames]{xcolor}
\usepackage{color}
\usepackage[usenames,dvipsnames]{xcolor}
\usepackage[colorlinks=true,linkcolor=Magenta,linktocpage=true,plainpages=false,allcolors=Magenta]{hyperref}
\usepackage{float}
\usepackage{epstopdf}
\usepackage{lipsum}
\usepackage{natbib}
\usepackage{doi}
\usepackage{soul}

\theoremstyle{definition}
\newtheorem{notation}{Notation}
\newtheorem*{proof*}{Proof}
\theoremstyle{plain}
\newtheorem{lemma}{Lemma}
\newtheorem{corollary}{Corollary}

\begin{document}

\setlength{\belowdisplayskip}{2pt} \setlength{\belowdisplayshortskip}{2pt}
\setlength{\abovedisplayskip}{2pt} \setlength{\abovedisplayshortskip}{2pt}

\begin{frontmatter}

\title{One and two sample Dvoretzky-Kiefer-Wolfowitz-Massart type inequalities for differing underlying distributions}

\author[1]{Nicolas G. Underwood\corref{cor1}}
\ead{mailto:NUnderwood@lincoln.ac.uk}
\author[1]{Fabien Paillusson} 
\ead{FPaillusson@lincoln.ac.uk}

\cortext[cor1]{Corresponding author}

\affiliation[1]{organization={School of Engineering and Physical Sciences, University of Lincoln},
            addressline={Brayford pool}, 
            city={Lincoln},
            postcode={LN6 7TS}, 
            country={United Kingdom}}

\begin{abstract}
Kolmogorov-Smirnov (KS) tests rely on the convergence to zero of the KS-distance $d(F_n,G)$ in the one sample case, and of $d(F_n,G_m)$ in the two sample case.
In each case the assumption (the null hypothesis) is that $F=G$, and so $d(F,G)=0$.
In this paper we extend the Dvoretzky-Kiefer-Wolfowitz-Massart inequality to also apply to cases where $F \neq G$, i.e.~when it is possible that $d(F,G) > 0$.
\end{abstract}

\end{frontmatter}

\section{Introduction}
Let $X_1,\dots,X_n$ and $Y_1,\dots,Y_m$ denote samples of univariate i.i.d.~random variables with continuous underlying cumulative distributions $F(x)$ and $G(x)$ respectively.
Denoting the indicator function $I(X_i\leq x)$, which returns 1 if $X_i \leq x$ and zero otherwise, 
the corresponding empirical distributions $F_n(x):=\frac{1}{n}\sum_{i=1}^nI(X_i\leq x)$ and $G_m(x):=\frac{1}{m}\sum_{i=1}^mI(Y_i\leq x)$ are step functions that approximate $F(x)$ and $G(x)$ for sample size $n$ and $m$ respectively.
The Kolmogorov-Smirnov (KS) distance, for instance $d(F_n,F):=\sup_x|F_n(x)-F(x)|$, is a sup-norm metric defined between two distributions, which takes values in $[0,1]$, and which quantifies the similarity of these distributions.
KS-type statistical hypothesis tests are based on a number of related theorems quantifying the convergence to zero with sample size of $d(F_n,G)$ in the one sample case, and of $d(F_n,G_m)$ in the two sample case for $F=G$ [and so $d(F,G)=0$]. 
In this paper, we find an analogous convergence to non-zero values of both $d(F_n,G)$ [Proposition \hyperlink{thm1}{1}], and $d(F_n,G_m)$ [Propositions \hyperlink{thm2a}{2a} and \hyperlink{thm2b}{2b}] in the case that $F\neq G$, so that $d(F,G) >0$.
 
\section{Convergence to zero}
One sample KS-tests were originally performed using a theorem by \cite{K33} on the asymptotic convergence of $d(F_n,F)$ to zero. Under the assumption that $F(x)$ is continuous, this states that
\begin{align}\label{kolmogorov}
\lim_{n\to\infty}\text{Prob}\left(\sqrt{n}d(F_n,F)\leq z \right)= L(z),
\end{align}
where $L(z)$ is given by
\begin{align}
L(z)=1-2\sum_{\nu=1}^\infty(-1)^{\nu-1}e^{-2\nu^2z^2}=\frac{\sqrt{2\pi}}{z}\sum_{\nu=1}^{\infty}e^{-(2\nu-1)^2\pi^2/8z^2}.
\end{align}
As usual, the corresponding $p$-value, $p\approx1-L[\sqrt{n}d(F_n,F)]$, represents the prior probability of obtaining a value at least as extreme as $d(F_n,F)$, given the null hypothesis that the data $X_1,\dots,X_n$ is drawn from distribution $F$.
However this is an estimate which only becomes exact in the large sample size limit, $n\to\infty$.
As sample sizes may not in practice be infinite, this one-sample test has largely been superseded by a test using a non-asymptotic theorem by \cite*{DKW56}, which was later refined by \cite{M90}.
The DKWM inequality states that 
\begin{align}\label{DKWM}
\text{Prob}\left[\sqrt{n}d(F_n,F) > z\right] \ \leq\  2e^{-2z^2},
\end{align}
and holds for all $z$ and $n$.
This provides a sharp upper bound to the $p$-value, $p\leq2 \exp\left[-2nd(F_n,F)^2\right]$.
Note that although this holds for all $z$, it may be said to have a "useful range" of $z\geq \sqrt{(\log 2)/2}$ or equivalently $d(F_n,F)\geq\sqrt{\log 2/(2n)}$, as outside of this range the probability upper bound is greater than unity.

Two-sample tests were traditionally performed using Smirnov's two-sample asymptotic convergence theorem \citep{S39,F49}.
This states that if $F$ and $G$ are the same continuous distribution, then
\begin{align}
\lim_{N\to\infty}\text{Prob}\left(\sqrt{N}d(F_n,G_m) \leq z \right)= L(z),
\end{align}
where $N:=nm/(n+m)$, and where the limit is understood to be taken such that the ratio $n/m$ tends to a constant.
The corresponding $p$-values, $p\approx1-L\left[\sqrt{N}d(F_n,G_m)\right]$, are prior probabilities of obtaining a value at least as extreme as $d(F_n,G_m)$, given the null hypothesis, continuous $F=G$.
In the same manner as Kolmogorov's one sample $p$-values, these are only approximate for finite sample sizes, becoming exact in the limit $N\rightarrow\infty$.
To account for finite sample sizes, a two-sample analogue to the DKWM inequality has been studied in some depth by \cite{WD12}.
Their results contain a number of detailed circumstances that we will not reexpress here, however for the case where $m=n$ so that $N=n/2$, it was found that
\begin{align}\label{WD}
\text{Prob}\left[\sqrt{\frac{n}{2}}d(F_n,G_n) > z\right] \ \leq\ Ce^{-2z^2},
\end{align}
where $C=2.16863$ for $n\geq 4$, and $C=2$ for $n\geq 458$. This bounds $p$-values as $p\leq C \exp\left[-nd(F_n,G_n)^2\right]$, and similarly to the DKWM inequality, has the useful range $z\geq \sqrt{(\log C)/2}$ or equivalently $d(F_n,G_n)\geq\sqrt{(\log C)/n}$.

\section{Convergence of $d(F_n,G)$ and $d(F_n,G_m)$ to finite values for $F\neq G$.}
Proposition \hyperlink{thm1}{1} adapts the DKWM inequality to place a non-asymptotic bound on the one-sample convergence of $d(F_n,G)$ to $d(F,G)$, for potentially distinct $F$ and $G$.
It follows from the metric property of the KS-distance, and possesses the same rate of convergence as $d(F_n,F)$ to zero as quantified by the DKWM inequality. 
\newtheorem*{thm1}{Proposition 1}
\begin{thm1}\hypertarget{thm1}{}
Let the empirical distribution $F_n$ and the cumulative distributions $F,G$ be defined as above. Then the KS-distance $d(F_n,G)$ converges to $d(F,G)$, with the rate of this convergence satisfying the inequality
\begin{align}\label{one_sample_inequality}
\mathrm{Prob}\left[\sqrt{n}|d(F_n,G)-d(F,G)| > z\right]\ \leq\ 2e^{-2z^2}.
\end{align}
\end{thm1}
\begin{proof*}
KS-distances are metric, and thus satisfy triangle inequalities, in particular the reverse triangle inequality,
\begin{align}
|d(F_n,G)-d(F,G)|\ \leq\ d(F_n,F).
\end{align}
For a given $y\in[0,1]$, as $|d(F_n,G)-d(F,G)|>y$ implies $d(F_n,F)>y$, it is evident that
\begin{align}
\mathrm{Prob}\big[|d(F_n,G)-d(F,G)| > y\big]\ \leq\ \mathrm{Prob}\left[d(F_n,F) > y\right].
\end{align}
Making the substitution $y=z/\sqrt{n}$, and using the DKWM inequality [Eq.~\eqref{DKWM}] on the right hand side, completes the proof. \qed
\end{proof*}
\hspace{-\parindent}Propositions \hyperlink{thm2a}{2a} and \hyperlink{thm2b}{2b} each place a non-asymptotic bound on the two-sample convergence of $d(F_n,G_m)$ to $d(F,G)$ given distinct $F$ and $G$. 
For equal sample size $m=n$, inequality \eqref{two_sample_inequality2_n=m} is a stronger bound than inequality \eqref{two_sample_inequality1_n=m} for $z>1.054$, corresponding to $p$-values of less than 78.4\%.
\newtheorem*{thm2a}{Proposition 2a}
\begin{thm2a}\hypertarget{thm2a}{}
Let distributions $F_n$, $G_m$, $F$, and $G$ be defined as above. The convergence of the KS-distance $d(F_n,G_m)$ to $d(F,G)$ satisfies the inequality
\begin{align}\label{two_sample_inequality1}
&\mathrm{Prob}\Big[|d(F_n,G_m)-d(F,G)| > z\Big]\nonumber\\
&\leq\frac{n}{m+n}2^{1-\frac{m}{n}}e^{-2mz\left(z-\sqrt{\frac{\log 4}{n}}\right)}
+\sqrt{32\pi}z\frac{mn}{(m+n)^{3/2}}e^{-\frac{2mnz^2}{m+n}}\mathrm{erf}\left(\sqrt{\frac{2m^2z^2}{m+n}}-\sqrt{\frac{(m+n)\log2}{n}}\right)+n\leftrightarrow m
\end{align}
for all $z\geq\sqrt{\log(2)/2}(n^{-1/2}+m^{-1/2})$.
In particular, for $n=m$ this may be expressed analogously to Eq.~\eqref{WD} as
\begin{align}\label{two_sample_inequality1_n=m}
\mathrm{Prob}\left[\sqrt{\frac{n}{2}}|d(F_n,G_n)-d(F,G)| > z\right]\ \leq\ \alpha_1(z):=
e^{-4z\left(z-\sqrt{\log 2}\right)}+\sqrt{32\pi}ze^{-2z^2}\mathrm{erf}\left[\sqrt{2}\left(z-\sqrt{\log 2}\right)\right]
\end{align}
for all $z\geq \sqrt{\log2}$.
\end{thm2a}

\begin{figure}[H]\label{functions}
\vspace{-10pt}
\begin{center}
\includegraphics[width=0.56\linewidth]{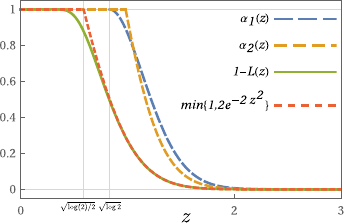}
\end{center}
\vspace{-10pt}
\caption{A comparison of functions $\alpha_1(z)$ and $\alpha_2(z)$ [defined in Propositions \protect\hyperlink{thm2a}{2a} and \protect\hyperlink{thm2b}{2b}] with $L(z)$ [used in Kolmogorov's and Smirnov's theorems] and $2e^{-2z^2}$ [used in the DKWM inequality, Wei and Dudley's inequality, and in Proposition \protect\hyperlink{thm1}{1}].}
\end{figure}

\newtheorem*{thm2b}{Proposition 2b}
\begin{thm2b}\hypertarget{thm2b}{}
With the same considerations as Proposition \hyperlink{thm2a}{2a}, the convergence satisfies
\begin{align}\label{two_sample_inequality2}
&\mathrm{Prob}\Big[|d(F_n,G_m)-d(F,G)| > z\Big]\nonumber\\
&\leq\min\left\{1,\left(1+\frac{n}{m+n}\right)2^{-\frac{m}{n}}e^{-2mz\left(z-\sqrt{\frac{\log4}{n}}\right)} 
+\sqrt{8\pi}z\frac{mn}{(m+n)^{3/2}}e^{-\frac{2mnz^2}{m+n}}\mathrm{erf}\left(\sqrt{\frac{2m^2z^2}{m+n}}-\sqrt{\frac{(m+n)\log2}{n}}\right)
+n\leftrightarrow m\right\}
\end{align}
for all $z\geq\sqrt{\log(2)/2}(n^{-1/2}+m^{-1/2})$.
For $n=m$, this may be expressed analogously to Eqs.~\eqref{WD} and \eqref{two_sample_inequality1_n=m} as
\begin{align}\label{two_sample_inequality2_n=m}
\mathrm{Prob}\left[\sqrt{\frac{n}{2}}|d(F_n,G_n)-d(F,G)| > z\right]\ \leq\ \alpha_2(z):=
\min\left\{1,\ \frac{3}{2}e^{-4z\left(z-\sqrt{\log 2}\right)}+\sqrt{8\pi}ze^{-2z^2}\mathrm{erf}\left[\sqrt{2}\left(z-\sqrt{\log 2}\right)\right]\right\}
\end{align}
for all $z\geq \sqrt{\log2}$. (Although this restriction is somewhat unnecessary as $\alpha_2=1$ for $z<0.983$.)
\end{thm2b}
\hspace{-\parindent}Figure \ref{functions} contrasts functions $\alpha_1(z)$ and $\alpha_2(z)$ [defined in Propositions \hyperlink{thm2a}{2a} and \hyperlink{thm2b}{2b}] with $L(z)$ [used in Kolmogorov's and Smirnov's theorems] and $2e^{-2z^2}$ [used in the DKWM and WD inequalities, and in Proposition \hyperlink{thm1}{1}].
\begin{lemma}\label{lemma1}
KS-distances satisfy the inequality,
\begin{align}\label{triangle}
|d(F_n,G_m)-d(F,G)|\ \leq\ d(F_n,F)+d(G_m,G).
\end{align}
\end{lemma}
\begin{proof*}
As KS-distances are metric, we may use triangle inequalities to find an upper bound on $d(F_n,G_m)$ as $d(F_n,G_m)\leq d(F_n,F)+d(F,G_m)\leq d(F_n,F)+d(G_m,G)+d(F,G)$.
A similar method applied to $d(F,G)$ results in the lower bound $d(F_n,G_m)\geq d(F,G)-d(F_n,F)-d(G_m,G)$.
Together these bounds may be recast as Eq.~\eqref{triangle}.\qed
\end{proof*}
\begin{corollary}\label{corollary1}
From Lemma \ref{lemma1}, it follows that for a given $z$, the statement $\big|d(F_n,G_m)-d(F,G)\big|>z$ implies the statement $d(F_n,F)+d(G_m,G)>z$. Hence,
\begin{align}\label{two_sample_prob_inequality}
\text{Prob}\Big[\big|d(F_n,G_m)-d(F,G)\big| > z \Big]\,\,\,\leq\,\,\, \text{Prob}\Big[d(F_n,F)+d(G_m,G) > z \Big].
\end{align}
\end{corollary}
\begin{notation}\label{notation1}
In the following, we abbreviate $x_F:=d(F_n,F)$ and $x_G:=d(G_m,G)$, and denote the (potentially unknown) distributions satisfied by these as $\rho_F(x_F)$ and $\rho_G(x_G)$. It is also useful to define the function 
\begin{align}\label{mu_1}
\mu_1(n,x) := \begin{cases}
8nxe^{-2nx^2}& \ \text{for}\ x\geq\sqrt{\frac{\log 2}{2n}}\\
0 & \mathrm{otherwise}
\end{cases},
\end{align}
which we use to  recast the DKWM inequality [Eq.~\eqref{DKWM}] in the form
\begin{align}\label{recast_DKWM}
\int_0^z\rho_F(x_F)\mathrm{d}x_F\geq \int_0^z\mu_1(n,x_F)\mathrm{d}x_F=\max\left(0,1-2e^{-2nz^2}\right).
\end{align}
\end{notation}
\begin{lemma}\label{lemma2}
For $z\leq 1$,
\begin{align}\label{integral}
\mathrm{Prob}\left[x_F+x_G\leq z\right]
\ \geq \ \int_{\omega_1} \mu_1(n,x_F)\mu_1(m,x_G)\mathrm{d}x_F\mathrm{d}x_G,
\end{align}
where $\omega_1:=\{(x_F,x_G):x_F+x_g\leq z\}$.
\end{lemma}
\begin{proof*}
As $x_F$ and $x_G$ are independent, the probability $\mathrm{Prob}\left[x_F+x_G\leq z\right]$ may be formally expressed as an integral of the combined probability distribution $\rho_F(x_F)\rho_G(x_G)$ over $\omega_1$. For $z\leq 1$, this is
\begin{align}
\mathrm{Prob}\left[x_F+x_G\leq z\right]=
\int_0^z\int_0^{z-x_G}\rho_F(x_F)\rho_G(x_G)\mathrm{d}x_F\mathrm{d}x_G\ \geq\ 
\int_0^z\int_0^{z-x_G}\mu_1(n,x_F)\rho_G(x_G)\mathrm{d}x_F\mathrm{d}x_G,
\end{align}
where we have applied the DKWM inequality [Eq.~\eqref{recast_DKWM}] to the $x_F$ integral. By switching the order of integration, we may once more apply the DKWM inequality, this time to the $x_G$ integral,
\begin{align}
\int_0^z\int_0^{z-x_F}\mu_1(n,x_F)\rho_G(x_G)\mathrm{d}x_G\mathrm{d}x_F\ \geq\ 
\int_0^z\int_0^{z-x_F}\mu_1(n,x_F)\mu_1(m,x_G)\mathrm{d}x_G\mathrm{d}x_F.
\end{align}
\qed
\end{proof*}
\theoremstyle{definition}
\newtheorem*{prf2a}{Proof of Proposition \hyperlink{thm2b}{2a}}
\begin{prf2a}
The proof of Proposition \hyperlink{thm2a}{2a} follows from Corollary \ref{corollary1} and Lemma \ref{lemma2}, and evaluation of integral \eqref{integral}.
We need only consider $z\leq 1$ as $|d(F_n,G_m)-d(F,G)|\leq 1$.
The restriction to $z\geq\sqrt{\log(2)/2}(n^{-1/2}+m^{-1/2})$ is a consequence of the piecewise definition of $\mu_1$ in Eq.~\eqref{mu_1}.
\qed
\end{prf2a}
We now turn to the proof of Proposition \hyperlink{thm2b}{2b}, which follows a similar approach to that of Proposition \hyperlink{thm2a}{2a}. 
Rather than relying on the inequality of Lemma \ref{lemma1}, we make use of the following stronger inequality.
\begin{lemma}\label{lemma3}
KS-distances $d(F_n,G_m)=\sup_x|F_n-G_m|$ and $d(F,G)=\sup_x|F-G|$ satisfy the inequality
\begin{align}\label{triangle-ish}
\big|d(F_n,G_m)-d(F,G)\big|
\ \leq \ 
\max \big[\sup_x(F_n-F)+\sup_x(G_m-G),\,\sup_x(F-F_n)+\sup_x(G-G_m)\big].
\end{align}
\end{lemma}
\begin{proof*}
Consider the following procedure.
\begin{align}
&\left(F_n(x)-G_m(x)\right) - \label{process1}
\left(F(x)-G(x)\right) = 
\left(F_n(x)-F(x)\right) + 
\left(G(x)-G_m(x)\right)\\
&\left(F_n(x)-G_m(x)\right) - \label{process2}
\left|F(x)-G(x)\right| \leq
\left(F_n(x)-F(x)\right) + 
\left(G(x)-G_m(x)\right)\\
&\left(F_n(x)-G_m(x)\right) - \label{process3}
\sup_x\left|F-G\right| \leq
\sup_x\left(F_n-F\right) + 
\sup_x\left(G-G_m\right)\\
&\sup_x\left(F_n-G_m\right) - \label{process4}
\sup_x\left|F-G\right| \leq
\sup_x\left(F_n-F\right) + 
\sup_x\left(G-G_m\right)
\end{align}
Eq.~\eqref{process1} applies pointwise. The absolute value applied in Eq.~\eqref{process2} serves to introduce the inequality, and the sup's introduced in Eq.~\eqref{process3} serve only to potentially loosen the inequality. Eq.~\eqref{process4} is reached as Eq.~\eqref{process3} is valid for all $x$.\footnote{For the sake of brevity, there is a slight abuse of notation in Eq.~\eqref{process3}; only the $\left(F_n(x)-G_m(x)\right)$ has $x$ dependence.}
By making the replacements $F_n\leftrightarrow G_m$, $F\leftrightarrow G$ on the LHS and $F_n\leftrightarrow F$, $G_m\leftrightarrow G$ on the RHS of Eq.~\eqref{process1}, an identical procedure may be followed to find
\begin{align}
&\sup_x\left(G_m-F_n\right) - \label{process5}
\sup_x\left|G-F\right| \leq
\sup_x\left(F-F_n\right) + 
\sup_x\left(G_m-G\right).
\end{align}
Then since $\sup_x\left|G-F\right|=\sup_x\left|F-G\right|$ and $\sup_x\left|F_n-G_m\right|=\max\left[\sup_x\left(F_n-G_m\right),\sup_x\left(G_m-F_n\right)\right]$, Eqs.~\eqref{process4} and \eqref{process5} may be combined to give
\begin{align}
&\sup_x\left|F_n-G_m\right| - \label{process6}
\sup_x\left|F-G\right|\ \leq \ 
\max\big[
\sup_x\left(F_n-F\right) + 
\sup_x\left(G-G_m\right),\ 
\sup_x\left(F-F_n\right) + 
\sup_x\left(G_m-G\right)\big].
\end{align}
Finally, this whole process may be repeated with the replacements $F_n\leftrightarrow F$ and $G_m\leftrightarrow G$, whereby the resulting equation may be combined with Eq.~\eqref{process6} to give Eq.~\eqref{triangle-ish}.\qed
\end{proof*}
\begin{notation}
In the following we abbreviate as $x_F^+:=\sup_x(F_n-F)$, $x_F^-:=\sup_x(F-F_n)$, $x_G^+:=\sup_x(G_m-G)$, and $x_G^-:=\sup_x(G-G_m)$.
Each of these satisfies the one-sided version of the DKWM inequality [see Thm.~1 in \cite{M90}], e.g. $\mathrm{Prob}(x_F^+>z)\leq \exp\left(-2nz^2\right)$ for all $z\geq\sqrt{\log2/(2n)}$. 
(Strictly speaking, we also require that $z\geq 1.0841 n^{-2/3}$, which we shall take as a given.) 
Similarly to our introduction of function $\mu_1$ in Notation \ref{notation1}, we define function $\mu_2(n,x):=4nx\exp\left(-2nx^2\right)$, which bounds the unknown distribution $\rho_{F^+}(x_F^+)$ satisfied by $x_F^+$ in the sense that we may reexpress the DKWM inequality 
\begin{align}\label{one_sided_DKWM}
\int_0^z\rho_{F^+}(x_F^+)\mathrm{d}x_F^+\geq\int_0^z\mu_2(n,x_F^+)\mathrm{d}x_F^+
\end{align}
for $z\geq\sqrt{\log2/(2n)}$.
\end{notation}
\begin{corollary}\label{corollary2}
From Lemma \ref{lemma3}, it follows that for a given $z$, the statement $\big|d(F_n,G_m)-d(F,G)\big|>z$ implies the statement $\max \big(x_F^+ + x_G^-,\,x_F^- + x_G^+\big)>z$. Hence,
\begin{align}
\mathrm{Prob}\Big[\big|d(F_n,G_m)-d(F,G)\big|>z\Big]
\ \leq \ 
\mathrm{Prob}\left[\max \big(x_F^+ + x_G^-,\,x_F^- + x_G^+\big)>z\right].
\end{align}
\end{corollary}
\begin{lemma}\label{alpha_3}
For $\sqrt{\log(2)/2}\left(n^{-1/2}+m^{-1/2}\right)\leq z \leq 1$,
\begin{align}\label{integral2}
\mathrm{Prob}\left[x_F^+ + x_G^-\leq z\right]
\ \geq \ 
\alpha_3:=
\int_{\omega_2} \mu_2(n,x_F^+)\mu_2(m,x_G^-)\mathrm{d}x_F^+\mathrm{d}x_G^-,
\end{align}
where $\omega_2:=\left\{(x_F^+,x_G^-): x_F^+ + x_G^- \leq z,\ x_F^+\leq z-\sqrt{\log 2/(2m)},\ x_G^-\leq z-\sqrt{\log 2/(2n)}\right\}$, and where we have defined $\alpha_3$ for later use.
(It is also the case that $\mathrm{Prob}\left[x_F^- + x_G^+\leq z\right]\ \geq \ \alpha_3$.)
\end{lemma}
\begin{proof*}
For $z\leq 1$, $\mathrm{Prob}\left[x_F^+ + x_G^-\leq z\right]=\int_{\omega_3}\rho_{F^+}(x_F^+)\rho_{G^-}(x_G^-)\mathrm{d}x_F^+\mathrm{d}x_G^-$, where $\omega_3:= \left\{(x_F^+,x_G^-): x_F^+ + x_G^- \leq z\right\}$.
\sloppy We may apply the one-sided DKWM inequality [Eq.~\eqref{one_sided_DKWM}] to the $x_F^+$ integral provided its upper boundary exceeds or is equal to $\sqrt{\log2/(2n)}$. We may ensure this by restricting to the integration domain to $\omega_4:=\left\{(x_F^+,x_G^-): x_F^+ + x_G^- \leq z,\ x_G^-\leq z-\sqrt{\log 2/(2n)}\right\}$, which only serves to loosen the resulting inequality, $\mathrm{Prob}\left[x_F^+ + x_G^-\leq z\right]\geq\int_{\omega_4}\mu_2(x_F^+)\rho_{G^-}(x_G^-)\mathrm{d}x_F^+\mathrm{d}x_G^-$ .
We may now switch the order of integration and apply the one-sided DKWM inequality to the $x_G^-$ integral provided its upper bound exceeds or is equal to $\sqrt{\log2/(2m)}$. This results in the integration domain $\omega_2$ defined above.
As $x_F^+\geq0$ and $x_G^-\geq0$, this domain requires $z\geq\sqrt{\log(2)/2}(n^{-1/2}+m^{-1/2})$.\qed
\end{proof*}
\theoremstyle{definition}
\newtheorem*{prf2b}{Proof of Proposition \hyperlink{thm2b}{2b}}
\begin{prf2b}
\sloppy In light of Corollary \ref{corollary2}, we look to place an upper bound on 
$\mathrm{Prob}\left[\max \big(x_F^+ + x_G^-,\,x_F^- + x_G^+\big)>z\right]
=\mathrm{Prob}\left[x_F^+ + x_G^- > z\ \bigcup\ x_F^- + x_G^+ > z \right]
=1-\mathrm{Prob}\left[x_F^+ + x_G^- \leq z\ \bigcap\ x_F^- + x_G^+ \leq z \right]$.
Thus we may equivalently look for a lower bound on $\mathrm{Prob}\left[x_F^+ + x_G^- \leq z\ \bigcap\ x_F^- + x_G^+ \leq z \right]$.
Due to Lemma \ref{alpha_3}, we have both $\mathrm{Prob}\left[x_F^+ + x_G^- \leq z\right]\geq \alpha_3$ and $\mathrm{Prob}\left[x_F^- + x_G^+ \leq z \right] \geq \alpha_3$.
(We need only consider $z\leq 1$ as $|d(F_n,G_m)-d(F,G)|\leq 1$.)
If $\alpha_3< 1/2$, there is the potential for $\mathrm{Prob}\left[x_F^+ + x_G^- \leq z\ \bigcap\ x_F^- + x_G^+ \leq z \right]=0$.
If however $\alpha_3> 1/2$, the smallest the intersectional probability may be is $2\alpha_3 -1$. 
Thus, $\mathrm{Prob}\left[x_F^+ + x_G^- \leq z\ \bigcap\ x_F^- + x_G^+ \leq z \right]\geq\max\left(0,\ 2\alpha_3 -1\right)$, and so $\mathrm{Prob}\left[\max \big(x_F^+ + x_G^-,\,x_F^- + x_G^+\big)>z\right]\leq 1-\max\left(0,\ 2\alpha_3 -1\right)=\min\left(1,2-2\alpha_3\right)$.
Eq.~\eqref{two_sample_inequality2} then follows from evaluation of $\alpha_3$ [Eq.~\eqref{integral2}].\qed
\end{prf2b}

\section{Conclusion}
KS-distances that are calculated from data [$d(F_n,G)$ in the one-sample case and $d(F_n,G_m)$ in the two-sample case] may be regarded as estimates of the underlying KS-distance, $d(F,G)$. 
The inequalities we have derived place probability bounds on the accuracy of these estimates. 
In the one-sample case, our Proposition \hyperlink{thm1}{1} reduces to the DKWM-inequality when $d(F,G)=0$.
When $d(F,G)=0$ in the two-sample case however, our Propositions \hyperlink{thm2a}{2a} and \hyperlink{thm2a}{2b} do not provide as tight a bound as that found by \cite{WD12}.
In addition to this, some preliminary numerical investigation suggests that the rate of convergence of $d(F_n,G_m)$ to $d(F,G)$ may quicken when $d(F,G)\neq0$.
For this reason, we speculate that it may be possible to find a tighter inequality in the two-sample case. 
Finally we note that our propositions in effect remove the $d(F,G)=0$ null hypothesis of standard Kolmogorov-Smirnov significance tests. 
To what extent this may inform such significance testing is an open question.

\section*{Acknowledgements}
This work was funded by Leverhulme Trust Research Project Grant RPG-2021-039.

\bibliographystyle{abbrvnat}
\bibliography{/home/nick/Research/ergodicity_citations}
\end{document}